# On the non-existence of skew-Hadamard difference sets in certain non-abelian groups

Vitor Araujo Garcia*

*Departamento de Ciências Exatas e da Terra*
*Faculdade de Ciência e Tecnologia*
*Universidade Federal de Mato Grosso, Cuiabá, MT, Brazil*
vitor.garcia@ufmt.br

**Abstract**

A skew-Hadamard difference set (SHDS) in a finite group $G$ is a classical combinatorial object with deep connections to design theory, coding theory, group theory, and the construction of Hadamard matrices. Even though the abelian case has been extensively studied — with strong structural constraints known, such as the necessity of $G$ being a $p$-group for some prime $p \equiv 3 \pmod 4$ — there are still some open questions regarding existence of SHDSs for the abelian case. The non-abelian case remains largely unexplored, despite the known existence of non-abelian SHDSs. In this paper, we establish new necessary conditions on the order and structure of a finite group $G$ that admits an SHDS. These results provide the first general structural restrictions for SHDSs in non-abelian groups. In particular, we prove that if a group $G$ is nilpotent and admits an SHDS, then $G$ is a $p$-group. Our method makes use of the structure of the rational group algebra, and completely avoids the use of the group characters.

## 1 Introduction

Let $G$ be a finite group of order $v$. A $(v, k, \lambda)$ *difference set* in $G$ is a subset $D \subseteq G$ of size $k$ such that each non-identity element of $G$ can be written as $d_1 d_2^{-1}$ (with $d_1 \neq d_2$ in $D$) in exactly $\lambda$ distinct ways.

Such a set $D$ is called a *skew-Hadamard difference set (SHDS)* if $G$ is the disjoint union of $D$, $D^{(-1)}$ and $\{1\}$, where $D^{(-1)}$ denotes the set of inverses of the elements of $D$. In this case, the parameters must satisfy $v = 4n-1$, $k = 2n-1$ and $\lambda = n-1$ for some integer $n$. Moreover, the identity

$$DD^{(-1)} = n + \lambda G,$$

* *ORCID:* 0000-0001-6150-5633

holds in the rational group algebra of $G$, where $D$, $D^{(-1)}$ and $G$ denote the formal sums of elements. Finally, if $D$ is an SHDS, then $D^{(-1)}$ is another SHDS with the same parameters, and $DD^{(1)} = D^{(-1)}D = n + \lambda G$.

Classical constructions of skew-Hadamard matrices are based on SHDSs, such as the Paley construction. Skew-Hadamard matrices are very important and have been widely studied in the literature, with several methods of construction; see, for example, [9]. For a comprehensive reference on difference sets, we refer the reader to [7].

It is well known that, if $G$ is abelian and admits an SHDS, then $G$ must be a $p$-group (see [6]). There has been significant effort to restrict which abelian $p$-groups could admit an SHDS (see [1], [2] and [6]), and it is conjectured that a finite abelian group $G$ admits an SHDS if and only if $G$ is elementary abelian. In [10], it was shown that if there exists an abelian but non-elementary abelian group admitting an SHDS, then there are infinitely many such groups, via a recursive construction.

For the non-abelian case, some constructions of SHDSs are known for specific $p$-groups (see [3], [4], [11]). In the non-abelian setting, although all known examples of SHDSs lie in $p$-groups, to the best of our knowledge, a general proof that this is always the case remains open.

Interestingly, in the preprint [13], the author states in its abstract that it is known that "if a group $G$ admits an SHDS, then $|G| = p^{2\alpha+1}$", without specifying that this result is only proven for abelian groups. In fact, in the main body of the preprint, all results are developed under the assumption that $G$ is abelian.

In the abelian setting, all known structural restrictions rely on the character table of $G$, a method that does not extend well to the non-abelian case. To overcome this limitation, we adopt a rational group algebra approach that completely avoids the use of characters.

In the present paper, we derive necessary conditions on the order and structure of finite groups $G$ that admit SHDSs. To the best of our knowledge, these results represent the first general structural constraints for SHDSs in non-abelian groups, addressing a gap in the literature and offering new tools for future construction efforts. In particular, we prove that if $G$ is a nilpotent group admitting an SHDS, then $G$ is a $p$-group. This generalizes a known result from the abelian case.

We use the following notation:

Let $\zeta_d$ denote a primitive $d$-th root of unity.

If $G$ is a group, we denote by $G' = [G, G]$ the commutator subgroup of $G$. We denote by $\mathbb{Q}[G]$ the rational group algebra over $G$ and, if $S \subseteq G$ we write $S$ in $\mathbb{Q}[G]$ to denote the sum of all elements of $S$ in $\mathbb{Q}[G]$. If $H \leq G$ we define:

$$\widehat{H} = \frac{1}{|H|}H = \frac{1}{|H|}\sum_{h\in H} h.$$

And finally, if $z = \sum_{g\in G} \alpha_g g \in \mathbb{Q}[G]$, we define $z^* = \sum_{g\in G} \alpha_g g^{-1}$ ($*$ is called the classical involution); it is clear that $(ab)^* = b^*a^*$.

### 1.1 Preliminaries

Here we present some algebraic tools that we need. If $R$ is any ring, we denote by $M_n(R)$ the ring of $n \times n$ matrices with coefficients in $R$. We begin with the following theorem:

**Theorem 1.1.** *[12, Theorem 3.4.9] Let $G$ be a finite group and $K$ a field such that* $\mathrm{char}(K) \nmid |G|$. *Then:*

1. *$K[G]$ is a direct sum of a finite number of two-sided ideals $\{B_i\}_{1\leq i\leq r}$, the simple components of $K[G]$;*
2. *Any two-sided ideal of $K[G]$ is a direct sum of some of the members of the family $\{B_i\}_{1\leq i\leq r}$;*
3. *Each simple component $B_i$ is isomorphic to a full matrix ring of the form $M_{n_i}(D_i)$, where $D_i$ is a division ring containing an isomorphic copy of $K$ in its center, and the isomorphism*

$$K[G] \overset{f}{\cong} \bigoplus_{i=1}^{r} M_{n_i}(D_i)$$

*is an isomorphism of $K$-algebras.*

The identity matrix $I_i$ in each component is a central idempotent in $B_i$. Let $e_i$ be the element in the group algebra $K[G]$ such that $f(e_i)$ is sent to the null matrix in each component $j \neq i$ and to $I_i$ in the $i^{th}$ component. Therefore:

1. $e_i e_j = 0$ if $i \neq j$ and $e_i^2 = e_i$ - we say that $e_i$ and $e_j$ are *orthogonal* for $i \neq j$;
2. $\sum_{i=1}^{r} e_i = 1$;
3. any two-sided ideal of $K[G]$ is generated by a finite sum of the elements $\{e_i\}_i$;
4. any central idempotent can be written as a sum of elements $e_i$.

The elements $e_i$ above are called the *primitive central idempotents* of $K[G]$, and any other central idempotent is a sum of certain $e_i$'s. In case $G$ is abelian, we have the following:

**Theorem 1.2.** *[12, Theorem 3.5.4] Let $G$ be a finite abelian group of order $v$, and let $K$ be a field such that* $\mathrm{char}(K) \nmid v$. *Then:*

$$K[G] \cong \bigoplus_{d|v} a_d K(\zeta_d),$$

*where $\zeta_d$ denotes a primitive root of unity of order $d$ and $a_d = \frac{n_d}{[K(\zeta_d):K]}$. In this formula, $n_d$ denotes the number of elements of order $d$ in $G$.*

The above theorems help us to understand the structure of $\mathbb{Q}[G]$. In particular, Theorem 1.2 suggests the possibility of connecting number theoretical facts about cyclotomic fields and the group algebra. In the next section, we relate the direct sum structures of $\mathbb{Q}[G]$ and $\mathbb{Q}[G/G']$ to obtain constraints on the order of a group $G$ that contains an SHDS.

*Remark* 1.3. To end this introductory section, we cite [5, Corollary 2.1]—in particular, this result tells us that, provided $G$ is a finite abelian group, then $\widehat{G}$ is a primitive central idempotent in the group algebra $\mathbb{Q}[G]$. If $x \in \mathbb{Q}[G]$, then $x\widehat{G} \in \mathbb{Q}$ is the sum of the coefficients of $x$, and we conclude that this particular idempotent corresponds to the summand $\mathbb{Q}$ in Theorem 1.2—if the order of $G$ is odd, there is only one such summand corresponding to $\mathbb{Q}$.

## 2 Constraints on the order of $G$

In this section we use algebraic and number theoretic tools to obtain restrictions on the order of groups $G$ that contain an SHDS. We begin with a key lemma, very similar to Theorem 1.2:

**Lemma 2.1.** *Let $G$ be a finite group. Then, as $\mathbb{Q}$-algebras,*

$$\mathbb{Q}[G](\widehat{G'} - \widehat{G}) \cong \bigoplus_{1 \neq d \mid [G:G']} a_d \mathbb{Q}(\zeta_d),$$

*where $\zeta_d$ is a $d$-primitive root of unity and $a_d = \frac{n_d}{[\mathbb{Q}(\zeta_d):\mathbb{Q}]}$, with $n_d$ the number of elements of $G/G'$ of order $d$.*

*Proof:* Let $T \subseteq G$ be a complete set of coset representatives of $G/G'$ and set $e_1 := (\widehat{G'} - \widehat{G})$. Note that $e_1$ is a central idempotent in the group algebra, and $e_1 = \widehat{G'}(1 - \widehat{G})$. Define:

$$\psi : \mathbb{Q}[G/G'] \to \mathbb{Q}[G]e_1$$
$$\sum_{t \in T} \alpha_t \bar{t} \mapsto \sum_{t \in T} \alpha_t t e_1.$$

It is easy to check that $\psi$ is well defined, because if $\overline{t_1} = \overline{t_2}$, then $t_1\widehat{G'} = t_2\widehat{G'}$. In addition, it is straightforward to verify that $\psi(ab) = \psi(a)\psi(b)$ and $\psi(a+b) = \psi(a) + \psi(b)$. Moreover, if $z = \sum_{g \in G} \alpha_g g e_1 \in \mathbb{Q}[G]e_1$, it is clear that

$$z = \sum_{t \in T} \left( \sum_{\bar{g} = \bar{t}} \alpha_g g \right) e_1 = \sum_{t \in T} \left( \sum_{\bar{g} = \bar{t}} \alpha_g \right) t e_1,$$

and $z$ is the image of $\sum_{t \in T} \left( \sum_{\bar{g} = \bar{t}} \alpha_g \right) \bar{t}$ under $\psi$, hence $\psi$ is a surjective homomorphism.

Let us determine its kernel. Let $A = \sum_{t\in T} \alpha_t \bar{t} \in \mathbb{Q}[G/G']$. We have

$$\begin{aligned}\psi\Big(\sum_{t\in T} \alpha_t \bar{t}\Big) &= \sum_{t\in T} \alpha_t t e_1\\ &= \sum_{t\in T} \alpha_t t \widehat{G'} - \sum_{t\in T} \alpha_t t \widehat{G}\\ &= \sum_{t\in T} \alpha_t t \frac{G'}{|G'|} - \sum_{t\in T} \alpha_t \frac{G}{|G|},\end{aligned}$$

since $tG = G$ for all $t \in G$.

Let us define the following elements in $\mathbb{Q}[G]$:

$$A_1 := \sum_{t\in T} \alpha_t t G'/|G'|,$$

$$A_2 := \sum_{t\in T} \alpha_t G/|G|,$$

so that $\psi(A) = A_1 - A_2$. We obtain that $A \in \ker(\psi)$ if and only if $A_1 = A_2$. Fix $g \in G$ and $t_g \in T$ such that $\overline{t_g} = \overline{g}$. We now study the coefficients of $g$ in $A_1$ and $A_2$, in order to compare them:

The coefficient of $g$ in $A_1$ is $\alpha_{t_g}/|G'|$, and the coefficient of $g$ in $A_2$ is $\sum_{t\in T} \alpha_t/|G|$. Therefore, $A \in \ker(\psi)$ if and only if $\alpha_{t_g}/|G'| = \sum_{t\in T} \alpha_t/|G|$ for all $g$, which is equivalent to

$$\alpha_{t_g} = \sum_{t\in T} \alpha_t |G'|/|G| \text{ , for all } g \in G.$$

Now, by comparing coefficients and noting that every $t \in T$ is a $t_g$ for some $g$, we get that $A \in \ker(\psi)$ if and only if $\alpha_t = \sum_{t'\in T} |G'|/|G| \alpha_{t'}$, for all $t \in T$, and this happens if and only if all $\alpha_t$'s are equal.

Thus, $\ker(\psi)$ is generated by the element $G/G'$ in the group algebra. Now, by Theorem 1.2 there is an isomorphism:

$$\mathbb{Q}[G/G'] \cong \bigoplus_{d | [G:G']} a_d \mathbb{Q}(\zeta_d),$$

where $\zeta_d$ is a $d$-primitive root of unity and $a_d = \frac{n_d}{[\mathbb{Q}(\zeta_d):\mathbb{Q}]}$, with $n_d$ the number of elements of $G/G'$ of order $d$.

Since $G/G'$ is abelian, by Remark 1.3, we find that the element $(G/G')/[G : G']$ is a primitive idempotent of $\mathbb{Q}[G/G']$. This element corresponds to the summand $\mathbb{Q}$ and generates the kernel of $\psi$.

Therefore:

$$\mathbb{Q}[G]e_1 \cong \mathbb{Q}[G/G']/\ker(\psi) \cong \bigoplus_{1\neq d|[G:G']} a_d \mathbb{Q}(\zeta_d).$$

□

**Theorem 2.2.** *Suppose $G$ is a finite group admitting an SHDS. Then $\sqrt{-|G|} \in \mathbb{Q}(\zeta_d)$ for all positive integers $d$ such that there is a cyclic subgroup of order $d$ in $G/G'$.*

*Proof:* Suppose $G$ is a group of order $4n-1$ and $D$ is an SHDS of $G$.

We define $x = D - D^{(-1)}$. Then:

$$\begin{aligned} x^2 &= D^2 - DD^{(-1)} - D^{(-1)}D + (D^{(-1)})^2 \\ &= D^2 + (D^{(-1)})^2 - 2DD^{(-1)}, \end{aligned}$$

since $DD^{(-1)} = D^{(-1)}D = (n + \lambda G)$, we have:

$$x^2 = D^2 + (D^{(-1)})^2 - 2(n + \lambda G).$$

Analogously, $(D + D^{(-1)})^2 = D^2 + (D^{(-1)})^2 + 2(n + \lambda G)$. But observe that $(D + D^{(-1)})^2 = (G-1)^2 = |G|G - 2G + 1$. From this, we obtain

$$D^2 + (D^{(-1)})^2 = (2n-1)G + 1 - 2n,$$

and substituting in the expression for $x^2$, and recalling $\lambda = n-1$ and $|G| = 4n-1$, this gives the key identity:

$$x^2 = G - |G|.$$

We now define the following elements in $\mathbb{Q}[G]$:

$$\begin{aligned} e_1 &= \widehat{G'}(1 - \widehat{G}), \\ e_2 &= 1 - e_1. \end{aligned}$$

We note that, if $G$ is abelian, $e_1 = (1 - \widehat{G})$ and $e_2 = \widehat{G}$.

It follows that $e_1$ and $e_2$ are idempotents (i.e., $e_i^2 = e_i$), orthogonal (i.e., $e_1 e_2 = 0$) and $e_1 + e_2 = 1$, since both $\widehat{G}$ and $\widehat{G'}$ are idempotents and $\widehat{G'}\widehat{G} = \widehat{G}$. And since $G'$ is a normal subgroup of $G$, we have that both $e_1$ and $e_2$ are central elements in $\mathbb{Q}[G]$ (i.e., they commute with all the elements in the group algebra). So, we can write

$$\mathbb{Q}[G] = \mathbb{Q}[G]e_1 \oplus \mathbb{Q}[G]e_2.$$

By Lemma 2.1, we have that the summand $\mathbb{Q}[G]e_1$ is isomorphic to a direct sum of fields $\mathbb{Q}(\zeta_d)$, for all $d \neq 1$ such that there is a cyclic subgroup of order $d$ in $G/G'$.

Moreover,

$$x^2 e_1 = (G - |G|)(\widehat{G'} - \widehat{G}) = (G - G) - |G|(\widehat{G'} - \widehat{G}) = -|G|e_1.$$

Therefore, $x^2$ projects onto the component corresponding to $e_1$ as $-|G|e_1$. In particular, $-|G|e_1$ must be a square in $\mathbb{Q}[G]e_1$.

Now, by Lemma 2.1, we have that $\mathbb{Q}[G]e_1$ is isomorphic to a direct sum of fields $\mathbb{Q}(\zeta_d)$, for all $d \neq 1$ such that there is a cyclic subgroup of order $d$ in $\mathbb{Q}[G/G']$, and this finishes our proof, since the element $1e_1 \in \mathbb{Q}[G]e_1$ corresponds, under this isomorphism, to the element $(1, 1, \ldots, 1)$ in the direct sum, and $x^2$ corresponds to $(-|G|, \ldots, -|G|)$ - thus $x$ corresponds to $(\pm\sqrt{-|G|}, \ldots, \pm\sqrt{-|G|})$. □

This theorem yields some non-existence results. We now prove some of them. First we recall a basic number theoretical fact: a quadratic field is an algebraic number field of degree two over $\mathbb{Q}$. If $l$ is an odd prime, then the only quadratic subfield of $\mathbb{Q}(\zeta_l)$ is either $\mathbb{Q}(\sqrt{l})$ or $\mathbb{Q}(\sqrt{-l})$; see [14, Section 6.5].

**Corollary 2.3.** *Let $G$ be a finite group. If any prime factor $\pi_1$ of $|G|$ has an odd multiplicity $\mu$, and $[G : G']$ has a prime factor $\pi_2 \neq \pi_1$, then $G$ does not admit an SHDS.*

*In particular, if $G$ admits an SHDS, then $G/G'$ is a $p$-group for some prime $p \equiv 3 \pmod 4$; furthermore, the multiplicity of $p$ in the primary decomposition of $|G|$ is odd.*

*Proof:* Let $\pi_1$ be a prime number dividing the order of $G$ with odd multiplicity, and $\pi_2$ another prime number dividing the order of $G/G'$. Suppose $G$ admits an SHDS.

Write $\mu = 2k + 1$, and we have:

$$\sqrt{-|G|} = \pi_1^k \sqrt{-\pi_1 M},$$

with $\gcd(M, \pi_1) = 1$.

By our theorem, we have $\pi_1^k\sqrt{-\pi_1 M} \in \mathbb{Q}(\zeta_{\pi_2})$, so $\sqrt{-\pi_1 M} \in \mathbb{Q}(\zeta_{\pi_2})$. This is impossible, since the only quadratic subfield of $\mathbb{Q}(\zeta_{\pi_2})$ is either $\mathbb{Q}(\sqrt{\pi_2})$ or $\mathbb{Q}(\sqrt{-\pi_2})$.

The last part of the corollary follows from the fact that $|G| \equiv 3 \pmod 4$. □

*Remark* 2.4. The result above generalizes a known fact (from [6]), that abelian groups admitting an SHDS are $p$-groups. In [6], the author needed the group to be abelian in order to be able to use nice properties of the group characters. In contrast, our method utilizes the abelianization $G/G'$ and avoids the use of group characters.

Now we write:

$$|G| = \prod p_i^{a_i} \prod q_j^{b_j},$$

with all $p_i$'s and $q_j$'s distinct primes and $p_i \equiv 3 \pmod 4$, $q_j \equiv 1 \pmod 4$ for all $i, j$.

The following results are immediate from Corollary 2.3.

**Corollary 2.5.** *With the notation above, if $b_j$ is odd for some $j$, then $G$ does not admit an SHDS.*

**Corollary 2.6.** *With the notation above, if the set* $\{i : a_i \text{ is odd}\}$ *has more than one element, then* $G$ *does not admit an SHDS.*

The following result is a direct consequence of Corollary 2.6:

**Corollary 2.7.** *Let* $G$ *be a group that admits an SHDS. If every prime factor of* $|G|$ *has odd multiplicity in the primary decomposition, then* $G$ *is a* $p$*-group for some* $p$.

The above results eliminate a broad class of groups from the possibility of admitting an SHDS. In fact, many additional corollaries could be derived from Theorem 2.2 and Corollary 2.3, covering further families of groups excluded by our necessary conditions.

We now present the case of finite nilpotent groups. This result not only generalizes the classical abelian case, but also shows that even in the broader class of nilpotent groups, the existence of an SHDS imposes strict structural constraints.

**Corollary 2.8.** *If* $G$ *is a finite nilpotent group admitting an SHDS, then* $G$ *is a* $p$*-group.*

*Proof:* We recall that if $G$ is a finite nilpotent group, then $G$ is the direct product $G = \prod_i P_i$, where each $P_i$ has prime power order. Hence $G/G' = \prod_i P_i/P_i'$, and since $|P_i|$ is odd (because $|G|$ is odd) we get that each factor is non-trivial (because the commutator subgroup of a $p$-group is always proper), and the result follows from Corollary 2.3. □

Our findings, together with known constructions in the literature, naturally lead to the following conjecture, which is now settled in the nilpotent case:

**Conjecture 2.9.** *If a finite group* $G$ *admits an SHDS, then* $G$ *is a* $p$*-group.*

It is worth noting that this conjecture remains open even for metacyclic groups. For example, our results do not apply to the group $C_{49} \rtimes C_3$ with trivial center, and we do not know whether it admits an SHDS.

If we suppose that $G$ is a $p$-group, trying to find exponent bounds similar to the bound given in [2] seems particularly challenging. Also, it is worth noting that in [4] the author constructed SHDSs for non-abelian groups of order $p^3$ and exponent $p > 3$, and according to [8], there is no such SHDS if $p = 3$, but there is one if we consider the non-abelian group $C_9 \rtimes C_3$ — so it looks a very hard problem to pose restrictions in the $p$-group case.